\documentclass[12pt]{amsart}
\usepackage{amssymb}
\usepackage{}
\usepackage{cases}

\usepackage{amsmath}
\usepackage{txfonts}
\theoremstyle{plain}
\newtheorem{Thm}{Theorem}

\newtheorem{Cor}[Thm]{Corollary}
\newtheorem{Pro}[Thm]{Proposition}
\newtheorem{Lem}[Thm]{Lemma}

\begin{document} 
\title[translating mean curvature flows]
{convexity and the Dirichlet problem of translating mean curvature flows}

\author{Li Ma}

\address{Ma Li: Zhongyuan Institute of mathematics and Department of mathematics \\
Henan Normal university \\
Xinxiang, 453007 \\
China} \email{lma@tsinghua.edu.cn}

\dedicatory{}
\date{May 26th, 2016}

\begin{abstract}

In this work, we propose a new evolving geometric flow (called translating mean curvature flow) for the translating solitons of hypersurfaces in $R^{n+1}$.
We study the basic properties, such as positivity preserving property, of the translating mean curvature flow. The Dirichlet problem for the graphical translating mean curvature flow is studied and the global existence of the flow and the convergence property are also considered.


\end{abstract}

\thanks{$^*$ The research is partially supported by the National Natural Science
Foundation of China (No. 11271111) and SRFDP 20090002110019. }
\maketitle

\section{Introduction}\label{sect1}
In this note, we propose a new evolving flow (called translating mean curvature flow) for the translating solitons of hypersurfaces in $R^{n+1}$. This flow is a modification of mean curvature flow with a translation by a fixed vector.
We study the basic properties of the translating mean curvature flow. The Dirichlet problem for the graphical translating mean curvature flow is studied and the global existence of the flow and the convergence property are also presented. This work can be considered as a continuation of our paper \cite{M}.

  We propose the translating mean curvature flow in the following way. Given a fixed nonzero vector $V\in R^{n+1}$. The translating mean curvature flow for translating soliton is defined as a one parameter family of properly immersed hypersurface $M_t=X(\Sigma,t)$, where $0<t<T$ and $X:\Sigma\times [0,T)\to  R^{n+1}$ evolved by the evolution equation
\begin{equation}\label{NF}
  X_t =\bar{H}(X)+V^N, \ \  t>0
\end{equation}
where $\bar{H}(X)$ is the mean curvature vector of the hypersurface $M_t$ at the position vector $X$ and $V^N$ is the normal component of the vector $V$. We denote $V^T=V-V^N$ the tangential part of the vector $V$. Recall that for the outer unit normal $\nu:=\nu(\cdot,t)$ on $M_t$, the mean curvature is defined by $H=div(\nu)$ and the mean curvature vector is $\bar{H}=-H\nu$. Let $(e_j)$ be a local orthonormal frame on $M_t$. Let
$a_{ij}=<D_{e{i}}e_j,\nu>$ and let $A=(a_{ij})$ be the second fundamental form on $M_t$. Then $H=-\sum a_{jj}$. Define
$\bar{X}=X-tV$. Then we have the mean curvature
$$
\bar{X}_t =\bar{H}(\bar{X}), \ \  t>0
$$
with the same initial hypersurface $X_0$. Therefore, many geometric properties such as convexity, mean convexity, are preserved by the flows. However, the global behaviors of two flows $X(t)$ and $\bar{X}(t)$ are different. Hence the flows (\ref{NF}) need to  be considered independently.

Applying the maximum principle (and Hamilton's tensor maximum principle) to derived evolution equations from (\ref{NF}) we obtain the following result.
\begin{Thm}\label{positive}
 Given a translating mean curvature flow $M_t$ with bounded second fundamental form $A$, $t\in [0,T)$ with $T>0$.

 (1). (i). If $<V,\nu>\geq 0$ on the initial hypersurface $M_0$, then $<V,\nu>\geq 0$ on the hypersurface $M_t$ for any $t>0$.
Similarly, assume that $H\geq 0$ on the initial hypersurface $M_0$. Then $H\geq 0$ on the hypersurface $M_t$ for any $t>0$.
 (ii). Assume that $A\geq 0$ on the initial hypersurface $M_0$. Then $A\geq 0$ on the hypersurface $M_t$ for any $t>0$.

 (2).  Assume that for some constant $\beta$,
$H-\beta<V,\nu>\leq 0$ on the initial hypersurface $M_0$ and $H-\beta<V,\nu><0$ at at least point in $M_0$.
Then
$H-\beta<V,\nu><0$ on $M_t$ for $t>0$.

(3). If we further assume $A+\beta\frac{<V,\nu>}{n}g\geq 0$ at the initial hypersurface $M_0$ and $A+\beta\frac{<V,\nu>}{n}g>0$ at at least one point $p\in M_0$, we have
$A+\beta\frac{<V,\nu>}{n}g>0$ on $M_t$ for $t>0$.
\end{Thm}

To derive this result, we shall do computations as in \cite{CM2}. As we have pointed out above, the property (1) can be derived from the mean curvature flow. For completeness, we give a full proof. Related Harnack inequalities for translating mean curvature flow similar to results in \cite{Ha} may be the same.

One example for hypersurfaces with $H-<V,\nu><0$ is
the graph of the parabolic function $u(x)=\frac{\lambda}{2} |x|^2$, where $x\in R^n$, $n\geq 2$ with $\lambda =1$ and $V=-e_{n+1}=(0,...0,-1)$. In this case, $Du(x)=x$, $v=\sqrt{1+|x|^2}$, $\nu=(-x,1)/v$,
$$
<\nu,V>=-1/v,
$$
and
$$
-H=div(\frac{x}{\sqrt{1+|x|^2}})=\frac{1}{v}(\frac{n+(n-1)|x|^2}{1+|x|^2})>\frac{1}{v}.
$$
One can compute that for $\lambda >0$ small we have $H-<V,\nu>>0$.

    The Dirichlet problem for translating solitons on convex domain has been studied by X.J.Wang (see Theorem 5.2 in \cite{W} from the viewpoint of Monge-Ampere equations. B.White \cite{Wh} has given a geometric measure theory argument for the existence of minimizers of the weighted area
    $$
    \int e^{-\lambda x_{n+1}}dA(x)
    $$
     amongst integral currents over the mean convex domain. Namely, letting $W$ be a bounded domain in $R^n$ with piecewise smooth mean convex boundary and letting $\Gamma$ be a smooth closed $(n-1)$ manifold in $\partial W\times R$ that is a graph-like. Then he has used the global defined radially symmetric solitons $y=\varphi(x)$ as barriers for the minimizing process of integral currents which lie in the region $R$ defined by
     $${R}=\{(x,y)\in \bar{W}\times R; b\leq y\leq \varphi(x)\}$$
     where $b=\inf \{y; (x,y)\in S\}$. We remark that his region
    $R$ (in the proof of Theorem 10 in \cite{Wh}) may be replaced by the region
    $$\check{R}=\{(x,y)\in \bar{W}\times R; \varphi(x)-C\leq y\leq \varphi(x)\}$$
    for suitable constant $C>0$. The choice of the lower barrier $\varphi(x)-C$ is nice in the sense that it is a sub-solution to the mean curvature soliton equation. One may get the minimizers by using BV functions. Our approach for the existence of translating solitons with the Dirichlet boundary condition on convex domains is the heat flow method. That is, we propose the translating mean curvature flow to get the solitons as the limits. The uniqueness and convexity of the translating solitons with convex boundary data $\phi$ remain as open questions.

    The Dirichlet problem for the graphical mean curvature flow on mean convex domains has been studied by G.Huisken \cite{Hu} and Lieberman \cite{L}. Their results show that the Dirichlet problem of the graphical mean curvature flow on mean convex domains has a global flow and it converges to a minimal surface at time infinity. Their result can not been directly applied to the following graphical translating mean curvature flow.
\begin{equation}\label{TMCF}
\partial_tu=\sqrt{1+|Du|^2}div (\frac{Du}{\sqrt{1+|Du|^2}})-1, \ \ \ in \ \ \Omega\times [0,\infty)
\end{equation}
with the Dirichlet boundary condition
$$
u=\phi, \ \ \ on \ \ \partial\Omega, t\geq 0
$$
and the initial condition
$$
u(x,0)=u_0(x), \ \ \ x\in \Omega.
$$
Here we assume $\Omega\subset R^n$ is a bounded domain with $C^2$ boundary, $\phi\in C^{2,\alpha}(\overline{\Omega})$, $u_0\in C^{2,\alpha}(\overline{\Omega})$, and $u_0=\phi$ on $\partial\Omega$. The flow (\ref{TMCF}) corresponds to the negative gradient flow of the weighted area functional
$$
F(u)=\int_{\Omega}\sqrt{1+|Du|^2}e^{u(x)}dx.
$$
If we let $f=u+t$, then $f$ satisfies
$$
\partial_tu=\sqrt{1+|Du|^2}div (\frac{Du}{\sqrt{1+|Du|^2}}), \ \ \ in \ \ \Omega\times [0,\infty)
$$
with the Dirichlet boundary condition
$$
f=\phi+t, \ \ \ on \ \ \partial\Omega, t\geq 0
$$
and the initial condition
$$
f(x,0)=u_0(x), \ \ \ x\in \Omega.
$$
Observe that the boundary condition now depends on time variable and the known result \cite{Hu} can not be applied directly to it.

We have the following result.
\begin{Thm}\label{Dirichlet} Assume $\Omega\subset R^n$ be a bounded convex domain with $C^2$ boundary. Assume that $\phi\in C^{2,\alpha}(\overline{\Omega})$, $u_0\in C^{2,\alpha}(\overline{\Omega})$, and $u_0=\phi$ on $\partial\Omega$.  Then the Dirichlet problem of (\ref{TMCF}) has a smooth solution and $u(\cdot,t)$ converges to the translating soliton with boundary data $\phi$ as $t\to\infty$.
\end{Thm}

The plan of this note is below. In section \ref{sect3} we discuss the positivity preserving properties of the general translating mean curvature flow. In section \ref{sect4} we consider the global existence of the Dirichlet problem of graphic mean curvature flows on bounded convex domains in $R^n$.

\section{positivity preserving property of the translating mean curvature flow}\label{sect3}

We shall use Hamilton's tensor maximum principle as below (see \cite{CL} for full statement and the proof).

\begin{Pro}\label{MP} Let $(M,g(t))$ be a one parameter family of complete noncompact Riemannian manifolds with bounded curvature.
Suppose $S = S_{ij}(x,t)dx^idx^j$ is a smooth time-dependent
symmetric 2-tensor field such that
$$
(\partial_t-\Delta_{g(t)})S\geq \nabla_XS+B(S,t)
$$
where $B(S,t)$ is locally Lipschitz in $(S,t)$ and $X=X(t)$ is a smooth time dependent vector field on $M$.
Assume that $B$ satisfies the null-eigenvector assumption in the sense that for some time-parallel vector field $v$ and at some point $x\in M$ such that if $S\geq 0$ and
$T(v,\cdot)=0$, then $B(S,t)(v,v)\geq 0$. Assume that $S\geq 0$ at the initial time $t=0$. Then $S\geq 0$ for all $t>0$.
\end{Pro}

Recall the following formulae for the flow $X_t=f\nu$ with local coordinates $(x_j)$ on $M_t$, we have for the evolving metric $g_{ij}=<\partial_{x_i}X,\partial_{x_j}X>$, outer unit normal $\nu$, and the second fundamental form $(a_{ij})$, we have
$$
\partial_tg_{ij}=-2fa_{ij},
$$
$$
\partial_t\nu=-\nabla f,
$$
and
$$
\partial_ta_{ij}=f_{ij}-fa_{ik}a_{kj}.
$$
We shall let $f=<V,\nu>-H$, which is our translating mean curvature flow case.

Let $(g^{ij})=(g_{ij})^{-1}$. As in \cite{CM2} and \cite{HP} we take $(e_i)$ to be the evolving frame on $M_t$ such that
$$
\partial_te_i=\frac{1}{2}g^{jk}\partial_tg_{ij}e_k=-fg^{jk}a_{ij}e_k.
$$
Then we have
$$
\partial_tg_{ij}=0.
$$
At a fixed point $p\in M_t$ we may assume that $<e_i,e_j>=\delta_{ij}$ and $\nabla_{e_i}e_j=0$. Then
$$
\partial_tA(e_i,e_j)=f_{ij}+fa_{ik}a_{kj}.
$$

Note that
$$
\nabla_i<V,\nu>=<V,D_{e_i}\nu>=-<V,e_k>a_{ik},
$$
and at $p$,
$$
\nabla_j\nabla_i<V,\nu>=-<V,e_k>a_{ik,j}-<V,D_{e_i}e_j>a_{ik}=a_{{ij},V^T}-<V,\nu>a_{ik}a_{kj}.
$$
Then we have
$$
\Delta <V,\nu>=\nabla_{V^T}H-<V,\nu>|A|^2.
$$
Since $\partial_t<V,\nu>=-<V,\nabla f>=-\nabla_{V^T}f$ and $f+H=<V,\nu>$, we get
$$
(\partial_t-\Delta )<V,\nu>=-\nabla_{V^T}f-\nabla_{V^T}H+<V,\nu>|A|^2
$$
$$
=-\nabla_{V^T}<V,\nu>+<V,\nu>|A|^2.
$$
That is,
\begin{equation}\label{Vp}
(\partial_t-\Delta )<V,\nu>=-\nabla_{V^T}<V,\nu>+<V,\nu>|A|^2.
\end{equation}

Recall the well-known formulae that
$$
(\Delta A)_{ij}=-|A|^2a_{ij}-Ha_{ik}a_{kj}-H_{ij}.
$$
Then we have
$$
(\partial_t A-\Delta A)(e_i,e_j)=(f+H)_{ij}+(f+H)a_{ik}a_{kj}+|A|^2a_{ij},
$$
which implies that for the normalized mean curvature flow,
$$
(\partial_t A-\Delta A)(e_i,e_j)=-\nabla_{V^T} a_{ij}+|A|^2a_{ij},
$$
that is,
\begin{equation}\label{second}
(\partial_t A-\Delta A)=-\nabla_{V^T} A+|A|^2A.
\end{equation}
Applying Hamilton's tensor maximum principle above (see also Proposition 12.31 in \cite{CL}) we know that $A\geq 0$ is preserved along the
translating mean curvature flow.
Note that by taking the trace of (\ref{second}), we have
$$
(\partial_t -\Delta)H=-\nabla_{V^T} H+|A|^2H.
$$
We can apply the scalar maximum principle to this equation and to (\ref{Vp}) too.
This gives the property (1) in Theorem \ref{positive}.

By these formulae for $A$, $H$, and $<V,\nu>$ we obtain that
$$
(\partial_t -\Delta) (A+\beta\frac{<V,\nu>}{n}g)=-\nabla_{V^T}(A+\beta\frac{<V,\nu>}{n}g)+|A|^2(A+\beta\frac{<V,\nu>}{n}g),
$$
and
$$
(\partial_t -\Delta) (H-\beta<V,\nu>)=-\nabla_{V^T}(H-\beta<V,\nu>)+|A|^2(H-\beta<V,\nu>).
$$

Define the operator
$$
L=\Delta-\nabla_{V^T}+|A|^2=\mathbf{L}+|A|^2.
$$
Then the above equations can be rewritten as
$$
(\partial_t -L)(A+\beta\frac{<V,\nu>}{n}g)=0
$$
and
$$
(\partial_t -L)(\beta<V,\nu>-H)=0.
$$
Applying the maximum principle (and Hamilton's tensor maximum principle as above) to above two equations we completes the proof of Theorem \ref{positive}.

One immediate consequence is the following pinching estimate.

\begin{Cor} Given a translating mean curvature flow $M_t$ with bounded second fundamental form $A$, $t\in [0,T)$ with $T>0$.
 Assume that for some uniform constants $\beta_1$ and $\beta_2$,
$\beta_1<V,\nu>\leq H\leq \beta<V,\nu>$ on the initial hypersurface $M_0$.
Then
$\beta_1<V,\nu>\leq H\leq \beta<V,\nu>$ on $M_t$ for $t>0$.
\end{Cor}

The proof is the same as (2) in Theorem \ref{positive}.

As in \cite{CM2} we have for any symmetric 2-tensor and $h$ a positive function on the manifold $M$,
$$
(\partial_t -\mathbf{L})|f|^2\leq 2<f,(\partial_t -\mathbf{L})f >,
$$
$$
(\partial_t -\mathbf{L})|\frac{f}{h}|^2\leq 2<\frac{f}{h},(\partial_t -\mathbf{L} \frac{f}{h}),
$$
and
$$
(\partial_t -\mathbf{L})\frac{f}{h}=\frac{(\partial_t -\mathbf{L})f }{h}- \frac{f(\partial_t -\mathbf{L})h}{h^2}+  \frac{2}{h}<\nabla h,\nabla \frac{f}{h}>.
$$
Then we have
$$
(\partial_t -\mathbf{L})|\frac{f}{h}|^2\leq 2<\nabla |\frac{f}{h}|^2,\nabla \log h>.
$$
Let, for some $\lambda$,
$$
B=\frac{A+\lambda \frac{<V,\nu>}{n}g}{\beta<V,\nu>-H}.
$$
\begin{Lem} Let $M_t\subset R^{n+1}$ be a one parameter family of hypersurfaces evolved by the translating mean curvature flow
(\ref{NF}). Assume that $\beta<V,\nu>-H>0$ on the initial hypersurface for some constant $\beta$, and $|A|^2$ are bounded on each $M_t$. Then
$$
(\partial_t -\mathbf{L})|B|^2\leq 2<\nabla |B|^2,\nabla \log (\beta<V,\nu>-H)>, \ \ on \ M_t.
$$
\end{Lem}

We now point out the geometric meaning of the operator $\Delta-\nabla_{V^T}+|A|^2$ on the hypersurface $M$.
Define the operator
$$
L=\Delta-\nabla_{V^T}+|A|^2,
$$
which is the Jacobian operator for the weighted volume
$$
F(M)=\int_M e^{-<V,X>}dX.
$$
Then, $F'=-\bar{H}-V^N=H\nu-V^N$.
In fact, for $X_t=X'=f\nu$ and $H'=\partial_t H$, we have
$$
H'=-\Delta f-|A|^2f
$$
and
$$
\nu'=-\nabla f.
$$
Then
$$
(H-<V,\nu>)'=-\Delta f-|A|^2f+<V,\nabla f>=-Lf.
$$
At the critical point of $F$ where
$$
H=<V,\nu>,
$$
we have
$$
F"=-\int_M <f,Lf>dm
$$
where $dm=e^{-<V,X>}dX$.

\section{The Dirichlet problem for the translating graphical mean curvature flow}\label{sect4}

Recall that
$\Omega\subset R^n$ is a bounded convex domain with $C^2$ boundary.

Note that the flow (\ref{TMCF}) corresponds to the negative gradient flow of the weighted area functional
$$
F(u)=\int_{\Omega}\sqrt{1+|Du|^2}e^{u(x)}dx.
$$
In fact,
$$
\delta F(u)\delta u=-\int_{\Omega}[div(\frac{Du}{\sqrt{1+|Du|^2}})-\frac{1}{v}]\delta ue^{u(x)}dx,
$$
where $v=\sqrt{1+|Du|^2}$. The functional $F(u)$ corresponds to the functional $F(M)$ with $V=-e_{n+1}$ in the previous section.

We point out a similarity between the translating mean curvature flow (\ref{TMCF}) and the graphical mean curvature flow. Fix any $t_0>0$.
Define
$U=u-t_0+t$.  Then $U$ satisfies the following
\begin{equation}\label{TMCF2}
\partial_tu=\sqrt{1+|Du|^2}div (\frac{Du}{\sqrt{1+|Du|^2}}), \ \ \ in \ \ \Omega\times [0,\infty)
\end{equation}
with the Dirichlet boundary condition
$$
U=\phi-t_0+t, \ \ \ on \ \ \partial\Omega, t\geq 0
$$
and the initial condition
$$
U(x,0)=u_0(x)-t_0, \ \ \ x\in \Omega.
$$
Define $Q=\sqrt{1+|DU|^2}$. Then $Q$ satisfies
\begin{equation}\label{grad}
\partial_tv=D_i(a^{ij}D_jv)+Hd^lD_lv-a^{ij}a^{kl}D_iD_kuD_jD_ku \cdot v, \ \ in \ \ \Omega\times (0,T)
\end{equation}
where
$a^i=Q^{-1}D_iU$,
$H=AU$,
and
$A^{ij}=\partial a^i/\partial p_j$. We shall use (\ref{grad}) to get the uniform gradient bound of $u$. We need to control
$\sup_{\partial\Omega} |Du|$ first.
Because of the equation (\ref{grad}) (being the same as the case of mean curvature flow) we believe the result of Theorem \ref{Dirichlet} should also be true for mean convex domains. However, we shall not discuss this in this note.

We now begin the proof of Theorem \ref{Dirichlet}.

\begin{proof} The existence of short time solution to (\ref{TMCF}) can be obtained by the standard method. Let $T>0$ be the maximal existence time of the solution $u(x,t)$. We claim that $T=+\infty$. To obtain this, we need to find a priori estimates for $\sup_\Omega |u|$ and $\sup_\Omega |Du|$.

Define
$$
Au=-div (\frac{Du}{\sqrt{1+|Du|^2}}).
$$
Let $w$ be the bowl soliton constructed by Altschuler-Wu \cite{AW}. Then we have
 $$
 -\sqrt{1+|Dw|^2}Aw=1, \ \ \ in \ \ \ R^n.
 $$

 Note that $w$ satisfies \ref{TMCF}. By adding to $w$
some uniform constant $C$ we may assume $w-C\leq -\sup_\Omega |u_0|$ and $w+C\geq \sup_\Omega |u_0|$. Using
$w\pm C$ as the barriers, we obtain that
$$
w-C\leq u\leq w+C, \ \ \ in \ \ \Omega\times [0,T).
$$
This gives us the uniform bound of $\sup_\Omega |u|$.

We now use the fact that the domain $\Omega$ is convex. Fix $p\in \partial\Omega$.
Recall that by the result of J.Serrin \cite{Se} or applying Cor. 14.3 in \cite{GT} to the operator
$$
Q(u)=-(\sqrt{1+|Du|^2})^3Au-(1+|Du|^2),
$$ we can construct barriers $\delta_+$ and $\delta_-$ such that $\delta_{\pm}(p)=\phi(p)$,
$$-(\sqrt{1+|D\delta_+|^2})^3A\delta_+\leq {1+|D\delta_+|^2}, \ \ \delta_+\geq \phi$$ and
$$ -(\sqrt{1+|D\delta_-|^2})^3A\delta_-\geq {1+|D\delta_-|^2}, \ \ \delta_-\leq \phi$$ in $\Omega$. We may also assume $\delta_-\leq u_0\leq \delta_+$ (see \cite{Hu}).

Applying the maximum principle to the evolution equation (\ref{TMCF}) we know that $\delta_-\leq u\leq \delta_+$ on $\Omega\times [0,T)$.
Hence we have at any time $t=t_0$ we have
$$
\delta_-(x)\leq u(x,t_0)\leq \delta_+(x)
$$
on $\Omega$. Since $p\in \partial\Omega$ is arbitrary, we know that there is a uniform constant $C_0$ depending only on $\partial\Omega$,$u_0$, and $\phi$ such that
$$
|Du|\leq C_0, \ \ \ on \ \ \partial\Omega\times [0,T).
$$
Applying the maximum principle the equation (\ref{grad}) for $Q$, we obtain the uniform bound for
$\sup_\Omega |Du|$. Once these are done, we then get the existence of the unique solution to the Dirichlet problem of (\ref{TMCF}) for all times $0<t<\infty$ with
the uniform gradient bound on $\sup_\Omega |Du|$. The standard parabolic equation theory \cite{LTU} guarantees uniform bounds of all higher derivatives of $u$.
Since $\partial_tu=0$ on $\partial\Omega$, by the equation we have  $H+\frac{1}{v}=0$ on $\partial\Omega$ and for $dm:=e^{u}dx$,
$$
\frac{d}{dt}\int_\Omega vdm=-\int_\Omega (H+\frac{1}{v})^2dm.
$$
Then
$$
\int_0^\infty\int_\Omega(H+\frac{1}{v})^2dm\leq \int_\Omega vdm(0).
$$
Using the uniform bound about $v$, we can conclude that
$\sup_\Omega |\partial_t u|$ and $\sup_\Omega |H-\frac{1}{v}|$ converges to zero uniformly as $t\to\infty$.
This completes the proof of Theorem \ref{Dirichlet}.

\end{proof}

\end{document}